\newtheorem{thm}{Theorem}
\newtheorem{lem}{Lemma}

\newtheorem{zam}{Remark}
\newtheorem{pred}{Proposition}

\newtheorem{gip}{Conjecture}
\newtheorem{sle}{Corollary}

\newcommand{\dvo}{{\it Proof}}

\documentclass[12pt]{article}
\usepackage{amsfonts}
\usepackage[all]{xy}
\title{On a Chisini Conjecture}\author{Vik.S. Kulikov
\thanks{Partly supported by RFFI
 (No. 96-01-00614) and INTAS (No. 96-0713).}}
\date{        }

\begin{document}
\maketitle
\begin{abstract}
Chisini's conjecture asserts that for a cuspidal curve $B\subset \mathbb P^2$ a generic morphism $f$ of a smooth projective surface onto $\mathbb P^2$ of degree $\geq 5$, branched along $B$, is unique up to isomorphism. We prove that if $\deg f$ is greater than the value of some function depending on the degree, genus, and number of cusps of $B$,  then the Chisini conjecture holds for $B$.  This inequality holds for many different generic morphisms. In particular, it holds for a generic morphism given by a linear subsystem of the $m$th canonical class for almost all surfaces with ample canonical class.
\end{abstract}   
\section*{Introduction}

Let $B\subset \mathbb P^2$ be an irreducible plane curve over $\mathbb C$ with ordinary cusps and nodes, as the only singularities. Denote by $2d$ the degree of $B$, and let $g$ be the genus of its desingularization, $c= \# \{ \mbox{cusps of}\, B\}$, and $n= \# \{ \mbox{nodes of} \, B\}$. We shall call $B$ {\it the discriminant curve} of a generic morphism if there exists a finite morphism $f:S\to \mathbb P^2$, $\deg f\geq 3$,  
satisfying the following conditions:

$(i)$ $S$ is a non-singular irreducible projective surface;

$(ii)$ $f$ is unramified over $\mathbb P^2 \setminus B$; 

$(iii)$  $f^{*}(B)=2R+C$, where $R$ is irreducible and non-singular, and $C$ is reduced; 

$(iv)$  $f_{\mid R}:R\to B$ coincides with the 
normalization of $B$. \newline 
We shall call such $f$ {\it a generic morphism}. 

Note that if $S\subset \mathbb P^r$, $f$ is the restriction to $S$ of a generic projection of $\mathbb P^r$ onto $\mathbb P^2$, and $B$ is the branch curve of $f$, then $(S,f)$  is a generic morphism and $B$ is its discriminant curve. 

Two generic morphisms $(S_1,f_1)$, $(S_2,f_2)$ with the same discriminant curve $B$ are said to be equivalent if there exists an isomorphism $\varphi : S_1 \to S_2$ such that 
$f_1=f_2\circ \varphi $. In the sequel, "$f$ is unique" means "$f$ is unique up to equivalence".

The following assertion is known as Chisini's Conjecture. 
\begin{gip}
Let $B$ be the discriminant curve of a generic morphism $f:S\to \mathbb P^2$ of degree $\deg f \geq 5$. Then, for $B$, the generic morphism $f$ is unique.
\end{gip}

If $B\subset \mathbb P^2$ is the dual curve of a smooth cubic, then $B$ is the discriminant curve of four generic non-equivalent morphisms (\cite{Chi}, \cite{Cat}). Three of them have degree four, and the last one has degree three.  This is the only known (up to now) example of a discriminant curve, for which there exist several non-equivalent generic morphisms.

In general case, as it follows from \cite{Cat}, the number of non-equivalent generic morphisms with a given discriminant curve $B$ is less than or equal to $2^{2g+c-1}$. 

B. Moishezon proved the Chisini Conjecture for 
the discriminant curves of generic projections of smooth 
hypersurfaces in $\mathbb P^3$. His proof is based on the presentation of the fundamental group of the complement in $\mathbb P^2$ of the discriminant curve of projection, obtained by him in \cite{Moi}. 
A short review of other results relating to the Chisini Conjecture, and of some attempts to prove it can be found in \cite{Cat}.

The main result of this paper is 

\begin{thm}
Let $B$ be the discriminant curve of a generic morphism
$f:S\to \mathbb P^2$ of $\deg f = N$. If 
\begin{equation}
N > \frac{4(3d+g-1)}{2(3d+g-1)-c}. \label{in}
\end{equation}
Then, for $B$, the generic morphism $f$ is unique and thus, the Chisini Conjecture holds for $B$.
\end{thm} 

Theorem 1 shows that if the degree of a generic morphism with given discriminant curve $B$ is sufficiently large, then this generic morphism is unique for $B$. Almost all generic morphisms interesting from algebraic geometric point of view satisfy this condition. More precisely, let $E=f^{*}(\mathbb P^1)$ be the preimage of a generic line $\mathbb P^1\subset \mathbb P^2$. In order to obtain the following theorems, which are consequences of the main result, we check inequality (\ref{in}) for morphisms given by three-dimensional subsystems of different linear systems $|E|$ on surfaces of different types.
\begin{thm}
Let $S$ be a surface of general type with ample canonical bundle $K_S$, $f:S\to \mathbb P^2$ a generic morphism such that $E\equiv mK_S$, $m\in \mathbb N$ ($\equiv $ means numerical equivalence). Then, for the discriminant curve $B$ of $f$, the generic morphism $f$ is unique except, possibly, for: 

1) $m=1,\, \, K_S^2=5,\, \, p_a=1,\, \, \deg B=20, \, \, g=51,\, c=108,\, \, n=12;$

2) $m=1,\, \, K_S^2=5,\, \, p_a=2,\, \, \deg B=20, \, \, g=51,\, c=96,\, \, n=24;$

3) $m=1,\, \, K_S^2=6,\, \, p_a=1,\, \, \deg B=24, \, \, g=61,\, c=132,\, \, n=60,$ \newline
where $p_a=\chi (\cal O_S)$ is the arithmetic genus of $S$.

In all the exceptional cases, if for $B$ there exist non-equivalent generic morphisms, then these morphisms have degree $\leq 6$.
\end{thm} 

\begin{thm}
Let $S$ be a Del Pezzo surface and $f:S\to \mathbb P^2$ be a generic morphism such that $E\in |-mK_S|$, $m\in \mathbb N$. Then, for the discriminant curve $B$ of $f$, the generic morphism $f$ is unique.
\end{thm} 

\begin{thm}
Let $f:S\to \mathbb P^2$ be any generic morphism of $S=\mathbb P^1\times \mathbb P^1$. Then, for the discriminant curve $B$ of $f$, the generic morphism $f$ is unique.
\end{thm} 

\begin{thm}
Let $S$ be a K3 surface and $f:S\to \mathbb P^2$ any generic morphism. Then, for the discriminant curve $B$ of $f$, the generic morphism $f$ is unique.
\end{thm}

\begin{thm}
Let $S$ be an Enriques surface and $f:S\to \mathbb P^2$ any generic morphism. Then, for the discriminant curve $B$ of $f$, the generic morphism $f$ is unique except, possibly, for $\deg f =4$. In the exceptional case, $\deg B=12, \, g=19,\, c=36,\, n=0$ and if such a morphism exists, then: \newline
1) for $B$, there exist at least two non-equivalent generic morphisms, \newline
2) any generic morphism $f'$ with such a discriminant curve $B$ has $\deg f'\leq 4$. \newline In particular, Chisini's Conjecture holds for the discriminant curves of the generic morphisms of Enriques surfaces.
\end{thm} 

\begin{thm}
Let $S$ be an abelian surface and $f:S\to \mathbb P^2$ any generic morphism. Then, for the discriminant curve $B$ of $f$, the generic morphism $f$ is unique except, possibly, for $\deg f =6$. In the exceptional case, $\deg B=18, \, g=28,\, c=72,\, n=36$ and if, for $B$, there exists a generic morphism $f'$ which is not equivalent to $f$, then $\deg f'\leq 6$.
\end{thm} 

\begin{thm}
Let $S\subset \mathbb P^N$ be a complete intersection and $f:S\to \mathbb P^2$ the restriction of a generic projection. Then, for the discriminant curve $B$ of $f$, the generic morphism $f$ is unique.
\end{thm} 

\begin{thm}
Let $S$ be a projective non-singular surface and $L$ an ample divisor on $S$, $f:S\to \mathbb P^2$ a generic morphism given by a three-dimensional subsystem $\{ E\} \subset |mL|$, $m\in \mathbb Q$, and $B$ its discriminant curve. Then there exists a constant $m_0$ (depending on $L^2, \, (K_S,L), \, K^2_S,\, p_a$) such that, for $B$, the generic morphism $f$ is unique if $m\geq m_0$. In particular, if $L=K_S$, then one can take $m_0=2$.
\end{thm} 

\begin{thm}
The Chisini Conjecture holds for the dual curve $B$ of a nodal plane curve except, possibly, for:

1) $\deg B= 30,\, \, g=10,\, \, c=72,\, \, n=324;$

2) $\deg B= 20,\, \, g=6,\, \, c=45,\, \, n=120;$

3) $\deg B= 18,\, g=5,\, \, c=39,\, \, n=92;$

4) $\deg B= 16,\, g=4,\, \, c=33,\, \, n=68.$

In all the exceptional cases, if, for $B$, there exist non-equivalent generic morphisms, then these morphisms have degree $\leq 6$.
\end{thm} 

\begin{thm}
The Chisini Conjecture holds for a curve $B$ of genus $g\leq 3$.
\end{thm} 

\begin{thm}
The Chisini Conjecture holds for $B$ satisfying the inequality 
$$ d> 3(g-1).$$
\end{thm} 

Unfortunately, in general case we have no a satisfactory description, purely in terms of algebraic geometry, of the set of discriminant curves with given degree, genus, and the number of cusps. But it is possible to give such a description in some particular cases. For example, in \cite{zar2}, Zariski showed that a sextic with 6 ordinary cusps is the discriminant curve of some generic morphism if and only if these 6 cusps lie on a conic, and he proved that there exist sextics with 6 cusps which do not lie on a conic. We also give, in terms of algebraic geometry, a description of discriminant curves of morphisms given by three-dimensional subsystems of the $m$th canonical class for $m\geq 21$ (see section 4.2 below). On the other hand, the set of the discriminant curves can be completely described in terms of the fundamental group of their complement in $\mathbb P^2$ (see Proposition 1 below).

In section 1, we recall some well-known facts on generic morphisms and their discriminant curves. Section 2 is devoted to the proof of Theorem 1. In section 3, we check inequality (\ref{in}) in different cases in order to prove Theorems 2 - 12. Section 4 contains a more detailed investigation of the case when a generic morphism is given by a linear subsystem of the $m$th canonical class. In section 5, we briefly discuss the question about the number of irreducible components of moduli space of discriminant curves with given degree, genus, and the number of cusps, and we apply Theorem 2 to find new examples of Zariski's pairs. 

This paper was written during my stay at the Max-Planck-Institut f\"{u}r Mathematik in Bonn. It is a pleasure to thank the Institut for its hospitality and financial support.

\section{Auxiliary results.}
{\bf 1.1.} Let $B$ be a discriminant curve of a generic morphism $f:S\to \mathbb P^2$, $\deg f=N$, $E=f^{*}(\mathbb P^1)$. We have $(E^2)_S=N$. 

\begin{lem}  $\deg B=2d$ is even and, consequently, $d\in \mathbb N$. 
\end{lem}
\dvo . By Hurwitz's formula, 
$$ 2g(E) -2= -2N +\deg B.$$
Hence $\deg B$ is even.

Since $g(E)\geq 0$, we have the following inequality
\begin{equation}
\deg f\leq d+1.
\end{equation}

Now we recall some inequalities contracting the genus, degree, and the number of cusps of $B$ and following from Pl\"{u}cker's formulas and Nori's result \cite{Nor}.  
For the dual curve $B^{*}$ of $B$, put $\delta =\deg B^{*}$, $\gamma = \# \{ \mbox{cusps of}\, \, B^{*}\}$ and $\nu = \# \{ \mbox{nodes of}\, \, B^{*}\}$. From Pl\"{u}cker's formulas
$$ 
\begin{array}{lll}
 \delta & = &  2d(2d-1)-2n-3c \, ;\\
 2d & = &  \delta (\delta -1)-2\nu -3\gamma \, ;  \\
 2g & = & (2d-1)(2d-2)-2n-2c \, ; \\
2g & = & (\delta -1)(\delta -2)-2\nu -2\gamma  
\end{array}
$$
it follows that 
$$ 
\begin{array}{lll}
 \delta & = &  4d-c +2g-2\, ;\\
 \gamma & = &  6d-2c+6g-6 . 
\end{array}
$$
Since $\delta \geq 0$ and $\gamma \geq 0$, we have
\begin{lem}  
$$
\begin{array}{lll}
c & \leq & 4d+2g-2\, ; \\
 c & \leq & 3d+3g-3\, .
\end{array}
$$
\end{lem}
\begin{sle}  $c< 2(3d+g-1) $.
\end{sle}
\begin{lem} $3d+g-1\leq 2c$. 
\end{lem}
\dvo . If $B$ is a discriminant curve of a generic morphism, then $\pi _1(\mathbb P^2\setminus B)$ is not abelian. Hence 
by \cite{Nor}, $B^2-6c\leq 2n$. For a plane cuspidal curve
$B$ of $\deg B=2d$ we have $g+c+n=(2d-1)(d-1)$. Hence 
$B^2-2n=4d^2-2n=2c+6d+2g-2\leq 6c$.

\begin{lem} \label{R} 
\begin{eqnarray}
R^2 & = & 2d^2-c-n = \label{R^21} \\
 & = & 3d+g-1   \label{R^2}
\end{eqnarray}
\end{lem}
\dvo . We have $K_S=-3E+R$. Hence (\ref{R^2}) follows from adjunction formula: 
$$\frac{(K_S+R,R)}{2}=\frac{(-3E+2R,R)}{2}=-3d+R^2=
g-1.$$
We obtain (\ref{R^21}) if instead of $g$ we substitute 
$g= (2d-1)(d-1)-c-n$. 

\begin{lem} \label{N<} 
\begin{equation}
\displaystyle N\leq \frac{4d^2}{3d+g-1},   \label{N<<}
\end{equation}
and (\ref{N<<}) is the equality if and only if either $E\equiv mK_S$ for some $m\in \mathbb Q^{*}$, or $K_S\equiv 0$.
\end{lem}
\dvo . By Hodge's Index Theorem,
$$
\left|
\begin{array}{cc}
E^2 & (E,R) \\
(E,R) &  R^2  
\end{array}
\right| 
=
\left|
\begin{array}{cc}
N & 2d \\
2d &  3d+g-1  
\end{array}
\right| 
\leq 0 ,
$$
i.e. 
$$\displaystyle N\leq \frac{4d^2}{3d+g-1} \, ,$$
and we have the equality if and only if $E$ and $R$ are linearly dependent in $NS(S)\otimes \mathbb Q$, where $NS(S)$ is the Neron-Severi group of $S$. In the last case, since $K_S=-3E+R$, either $E\equiv mK_S$ for some $m$, or $K_S\equiv 0$.

\begin{lem}  $$
\begin{array}{lll}
K_S^2 & = & 9N+2d^2-12d-c-n = \\
           & = & 9N-9d+g-1.
\end{array}
$$
\end{lem}
\dvo . $K_S^2 =(-3E+R,-3E+R)=9E^2-6(E,R)+R^2=9N-12d+R^2$.

\begin{lem} The topological Euler characteristic 
$$
\begin{array}{lll}
e(S) & = & 3N+4d^2-6d-3c-2n = \\
           & = & 3N+2(g-1)-c.
\end{array}
$$
\end{lem}
\dvo .
Using a generic pencil of lines in $\mathbb P^2$ and its preimage $\{ E_t \}$ in $S$, we get a formula for $e(S)$:
$$
\begin{array}{lll}
e(S)+N & = & \delta + 2(2-2g(E))= 2d(2d-1)-3c-2n -2(K_S+E,E)= \\
 & = & 2d(2d-1)-3c-2n-2(-2E+R,E)=
4N+4d^2-6d-3c-2n.
\end{array}
$$

{}From Noether's formula $K_S+e(S)=12p_a$ it follows
\begin{lem} The Euler characteristic of $\cal O_S$
$$
\begin{array}{lll}
\displaystyle p_a=1-q+p_g & = & N+\displaystyle     \frac{d(d-3)}{2}-\frac{c}{3}-\frac{n}{4}= \vspace{0.3cm} \\ 
& = & N+\displaystyle \frac{3g-3-9d-c}{12} .
\end{array}
$$
\end{lem}

\begin{sle} (\cite{Moi2}). 
$$
\begin{array}{lll}
 c & \equiv &  0 \,\, (\mbox{mod}\, 3)\, ;\\
 n & \equiv & 0 \, \, (\mbox{mod}\, 4)\, .
\end{array}
$$
\end{sle}
\begin{lem} The divisor $R$ is ample on $S$.
\end{lem}
\dvo . It suffices to show that $(R,\Gamma )>0$ for any irreducible curve $\Gamma $. If 
$(R,\Gamma )\leq 0$, then by Hodge's Index Theorem, $\Gamma ^2<0$ and $(R,\Gamma )=0$, since by Lemma \ref{R}, $R^2>0$ and $R$ is irreducible. If $\Gamma ^2<0$ and $\Gamma $ is irreducible, then $(K_S,\Gamma )\geq -1$, i.e. 
$-3(E,\Gamma )+(R,\Gamma )\geq -1$, which contradicts $(R,\Gamma )=0$.
\newline 
{\bf 1.2.} 
Let us fix $p\in \mathbb P^2\setminus B$ and denote by $\pi _1=\pi _1(\mathbb P^2\setminus B,\, p)$ the fundamental group of the complement of $B$. Choose any point $x\in B\setminus Sing\, B$ and consider a line $\Pi =\mathbb P^1\subset \mathbb P^2$ intersecting $B$ transversely at $x$. Let $\gamma \subset \Pi$ be a circle of small radius with center at $x$. If we choose an orientation on $\mathbb P^2$, then it defines an orientation on $\gamma $. Let $\Gamma $ be a loop consisting of a path $L$ in $\mathbb P^2\setminus B$ joining the point $p$ with a point $q\in \gamma$, the circuit in positive direction along $\gamma $ beginning and ending at $q$, and  a return to $p$ along the path $L$ in the opposite direction. Such loops $\Gamma$ (and the corresponding elements in $\pi _1$) will be called {\it geometric generators}. It is well-known that $\pi _1$ is generated by geometric generators, and any two geometric generators are conjugated in $\pi _1$, since $B$ is irreducible.

For each singular point $s_i$ of $B$ we choose a small neighborhood $U_i\subset \mathbb P^2$ such that $B\cap U_i$ is defined (in local coordinates in $U_i$) by equation $y^2=x^3$ if $s_i$ is a cusp, and $y^2=x^2$ if $s_i$ is a node. Let $p_i$ be a point in $U_i\setminus B$. It is well-known that if $s_i$ is a cusp, then $\pi _1(U_i\setminus B,p_i)$ is isomorphic to the braid group $\mbox{Br}_3$ of 3-string braids and generated by two geometric generators (say $a$ and $b$) satisfying the following relation 
$$aba=bab.$$ 
If $s_i$ is a node, then $\pi _1(U_i\setminus B,p_i)$ is isomorphic to $\mathbb Z\oplus \mathbb Z$ generated by two commuting geometric generators. 

Let us choose smooth paths $\gamma _i$ in $\mathbb P^2\setminus B$ joining $p_i$ and $p$. This choice 
defines homomorphisms $\psi _i:\pi _1(U_i\setminus B,p_i)\to \pi _1$. Denote the image  $\psi _i(\pi _1(U_i\setminus B,p_i))$ by $G_i$ if $s_i$ is a cusp, and $\Gamma _i$ 
if $s_i$ is a node.

A generic morphism of degree $N$ determines a homomorphism $\varphi :\pi _1 \to \mathfrak S_N$, where $\mathfrak S_N$ is the symmetric group. This homomorphism $\varphi $ is determined uniquely up to inner automorphism of $\mathfrak S_N$.
\begin{pred}
The set of the non-equivalent generic morphisms of degree $N$ possessing  the same  discriminant curve $B$ is in one to one correspondence with 
the set of the epimorphisms $\varphi :\pi _1(\mathbb P^2\setminus B)\to \mathfrak S_N$ (up to inner automorphisms of $\mathfrak S_N$) satisfying the following conditions:

$(i)$ for a geometric generator $\gamma $ the image $\varphi (\gamma )$ is a transposition in $\mathfrak S_N$; 

$(ii)$ for each cusp $s_i$ the image $\varphi (G_i)$ is  isomorphic to $\mathfrak S_3$ generated by two transpositions;

$(iii)$ for each node $s_i$ the image $\varphi (\Gamma _i)$
is isomorphic to $\mathfrak{S}_2 \times \mathfrak{S}_2$ generated by two commuting transpositions.
\end{pred}
\dvo . It is well-known that each homomorphism $\varphi :\pi _1 \to \mathfrak S_N$ defines a finite morphism $f:S\to \mathbb P^2$ of degree $N$, unramified over the complement of $B$, and such that $S$ is a normal surface, and vice versa. 

Condition $(i)$ is equivalent to one that $f^{*}(B)=2R+C$, where $R$ is irreducible and $C$ is reduced (cf. \cite{Kul1} and \cite{Kul}).

Conditions $(ii)$ and $(iii)$ are equivalent to one that $S$ and $R$ are non-singular at the points in $f^{-1}(s_i)$ (cf. \cite{Cat}). 

Since $S$ is irreducible and $\deg f =N$, $\varphi (\pi _1)$ must act transitively on the set $\overline{N}=\{1,...,\, N\}$, and since $\pi _1$ is generated by the geometric generators, $\varphi (\pi _1)$ is generated by some subset of transpositions. It is easy to check that subgroup of $\mathfrak S_N$, generated by some subset of transpositions acting on $\overline{N}$  transitively, must coincide with $\mathfrak S_N$. Therefore $\varphi $ must be an epimorphism.
\begin{zam}
By Proposition 1, if we allow $N=2$, then a two-sheeted covering $f:S\to \mathbb P^2$, branched along a non-singular curve $B\subset \mathbb P^2$, can be also considered as a generic morphism.
\end{zam}
\begin{zam}
If $B\subset \mathbb P^2$ is the dual curve of a smooth cubic, then by Proposition 1, using the presentation of $\pi _1(\mathbb P^2\setminus B)$ obtained in \cite{Zar}, it is easy to show that $B$ is the discriminant curve of exactly four generic non-equivalent morphisms.
\end{zam}
{\bf 1.3.} 
Let $s_j \in B$ be a cusp and $(x,y)$ local coordinates around $s_j$ such that $B\cap U_j$ is given by $y^2=x^3$. Choose a neighborhood  $V_j\subset S$ of $p_j\in R$, $f(p_j)=s_j$, such that $f=f_{\mid V_j}: V_j\to U_j$ is a three-sheeted covering ramified along $R\cap V_j$. It is well known that $f :V_j\to U_j=U$ is unique up to equivalence and, in particular, $f$ is equivalent to the standard covering $\overline{f} :\overline{V} \to U$ given by the normalized equation of third degree:
$$
\left\{
\begin{array}{l}
\overline{V}=\{ (w,x,y)\, \, \mid \, \, (x,y)\in U, \, \, w^3-3xw+2y=0\, \, \, \} \, ;\\
\overline{f}(w,x,y)=(x,y)\, .
\end{array}
\right.
$$
$\overline{V}$ is non-singular and $(x,w)$ are local coordinates in $\overline{V}$. The ramification divisor $\overline{R}=\{ (x,w)\in \overline{V}\, \, \mid \, \, x-w^2=0\, \} $ of $\overline{f}$ is smooth  and $\overline{f}^{-1}(B)=2\overline{R}+\overline{C}$, where $\overline{C}=\{ (x,w)\in \overline{V}\, \, \mid \, \, 4x-w^2=0\, \} $. Note that $\overline{R}$ is tangent to $\overline{C}$ at the origin $o=(0,0)$, and the intersection multiplicity of $\overline{R}$ and $\overline{C}$ at $o$ is equal to 2. 

It was mentioned above that $\pi _1(U\setminus B) \simeq \mbox{Br}_3=<\, a,b\, \, \mid \, \, aba=bab\, >$. Then $\overline{f} :\overline{V} \to U$ corresponds to the homomorphism $\overline{\varphi }:\pi _1(U\setminus B) \to \mathfrak S_3$ given by $\overline{\varphi }(a)=(1,2)$ and $\overline{\varphi }(b)=(2,3)$.

Put $\overline{W}=\{ (w_1,w_2,w_3)\, \, \mid \, \, w_1+w_2+w_3=0\, \} $. We have morphisms $\widetilde f : \overline{W}\to U$ and  $g : \overline{W}\to \overline{V}$ given by 
$$
\begin{array}{lll}
x & = & \displaystyle -\frac{1}{3}(w_1w_2+w_1w_3+w_2w_3) ,\vspace{0.3cm} \\
y & = & \displaystyle -\frac{1}{2}w_1w_2w_3 ,\\
w & = & w_1 
\end{array}
$$
such that $ \widetilde f =\overline{f}\circ g$, $\deg \widetilde f =6$ and $\deg g=2$. The morphism $g$ is a two-sheeted covering branched along $\overline{C}$, $g^{*}(\overline{C})=2\widetilde C_2$, where $\widetilde C_2$ is given in coordinates $(w_1,w_2)$ by $w_1+2w_2=0$; and  $g^{*}(\overline{R})=\widetilde R+\widetilde C_1$, where $\widetilde R$ and $\widetilde C_1$ are given by $w_1=w_2$ and $2w_1+w_2=0$ respectively.

Note that $\widetilde f$ corresponds to the homomorphism 
$\widetilde{\varphi }:\pi _1(U\setminus B) \to \mathfrak S_6=\mathfrak S(\mathfrak S_3)$ defined by $\overline{\varphi }$.

\section{Proof of Theorem 1.}
{\bf 2.1.} Assume that there exist two non-equivalent generic morphisms $(S_1,f_1)$ and $(S_2,f_2)$ with the same discriminant curve $B$, $\deg f_1=N_1$ and $\deg f_2=N_2$. Put $f^{*}_1(B)=2R_1+C_1$ and $f^{*}_2(B)=2R_2+C_2$.
Let us consider $$S_1\times _{\mathbb P^2}S_2=\{ \, (x,y)\in S_1\times S_2\, \,  \mid \, \, f_1(x)=f_2(y)\, \, \} $$
and let $X=\widetilde{S_1\times _{\mathbb P^2}S_2}$ be the normalization of $S_1\times _{\mathbb P^2}S_2$. Denote by $g_{1}:X\to S_1$, $g_{2}:X\to S_2$, and $f_{1,2}:X\to \mathbb P^2$ the corresponding natural morphisms. We have $\deg g_1=N_2$, $\deg g_2=N_1$, and  $\deg f_{1,2}=N_1N_2$.

\begin{pred} If $(S_1,f_1)$ and $(S_2,f_2)$ are non-equivalent, then $X$ is irreducible.
\end{pred}

\dvo . The morphism $f_{1,2}$ corresponds to the homomorphism $$\varphi _{1,2}=\varphi _1\times \varphi _2 : \pi _1 \to \mathfrak S_{N_1}\times \mathfrak S_{N_2}\subset \mathfrak S_{N_1N_2},$$ 
where $\varphi _1 : \pi _1 \to \mathfrak S_{N_1}$ (resp. $\varphi _2$) is an epimorphism corresponding to $f_1$ (resp. $f_2$). 
Put 
$G=\varphi _{1,2}(\pi _1)$. The group $G$, as a subgroup of $\mathfrak S_{N_1}\times \mathfrak S_{N_2}$, acts on $\overline{N_1}\times \overline{N_2}$. Without loss of generality we can assume that for some geometric generator $\gamma $ the image $\varphi _{1,2}(\gamma )=((1,\, 2),(1,\, 2))$ is a product of transpositions $(1,\, 2) \in\mathfrak S_{N_i}$.  Let $p_i: G\to \mathfrak S_{N_i}$ be the restriction to $G$ of the projection  $pr_i :\mathfrak S_{N_1}\times \mathfrak S_{N_2}\to \mathfrak S_{N_i}$. 

\begin{lem} Let $G$ be a subgroup of $\mathfrak S_{N_1}\times \mathfrak S_{N_2}$, $N_i>2$, such that $p_i: G\to \mathfrak S_{N_i}$ is an epimorphism for $i=1,\, 2$, and such that $((1,\,2), (1,\,2)) \in G$. Let $St_{(1,1)}\subset G$ be the stabilizer of  $(1,\,1)\in  \overline{N_1}\times \overline{N_2}$. Then the index of $St_{(1,1)}$ in $G$
$$( G:St_{(1,1)})=N_1N_2$$
except for the case when $N_1=N_2=N$ and $G= \Delta \subset \mathfrak S_{N}\times \mathfrak S_{N}$ (up to inner automorphism of one of factors), where $\Delta $ is the diagonal subgroup.
\end{lem}

\dvo . The inclusion $\overline{N_{i}-1}\simeq \{ 2,\, ...\, , \, N_i\}\subset \{ 1,\, 2,\, ...\, , \, N_i\}$ defines the embedding $\mathfrak S_{N_i-1}\subset \mathfrak S_{N_i}$. Then 
$$St_{(1,1)}=G\cap (\mathfrak S_{N_1-1}\times \mathfrak S_{N_2-1}).$$ Put 
$$H_1\times \{ e_2 \}= G\cap (\mathfrak S_{N_1}\times \{ e_2 \})\hspace{1cm} \mbox{and} \hspace{1cm} \{ e_1 \} \times H_2= G\cap (\{ e_1 \} \times \mathfrak S_{N_2}),$$
 where $e_i$ is the unit element of $\mathfrak S_{N_i}$. We note that $\ker p_2=H_1\times \{ e_2 \}$ and $\ker p_1=\{ e_1 \} \times H_2$. Therefore $H_1\times \{ e_2 \}$ and $\{ e_1 \} \times H_2$ are normal subgroups of $G$. Since $p_i: G\to \mathfrak S_{N_i}$ is an epimorphism for $i=1,\, 2$, $H_i$ is a normal subgroup of $\mathfrak S_{N_i}$. 

It is well known that  if $H$ is a normal subgroup of $\mathfrak S_{N}$, then either $H=\mathfrak S_{N}$, or $H=\mathfrak A_{N}$ is the alternating group, or $H= \{ e \}$, and if $N=4$, then there exist one more possibility: 
$H$ is the Klein four group 
$$K_4=\{ e,\, \,  (1,2)(3,4),\, \, (1,3)(2,4),\, \, (1,4)(2,3)\} .$$
Consider all possible cases.

{\it Case I :}  $H_1=\mathfrak S_{N_1}$. Since $p_2$ is an epimorphism, $G=\mathfrak S_{N_1}\times \mathfrak S_{N_2}$.
Therefore $\mid G\mid =N_1!N_2!$ and $\mid St_{(1,1)}\mid =(N_1-1)!(N_2-1)!$. Hence $( G:St_{(1,1)})=N_1N_2$.

{\it Case II :}  $H_1=\mathfrak A_{N_1}$. Since $p_2$ is an epimorphism and $\ker p_2=H_1\times \{ e \}$, 
$$\mid G\mid =\frac{N_1!N_2!}{2}.$$ 
Similarly, if we consider $p_1$, then we obtain $|\ker p_1|=N_2!/2$, hence $H_2=\mathfrak A_{N_2}$. Therefore $(\sigma _1, \sigma _2)\in G$ if and only if $\sigma _1 $ and $\sigma _2$ have the same sign. Hence $$\mid St_{(1,1)}\mid = \frac{(N_1-1)!(N_2-1)!}{2}$$ 
and $(G:St_{(1,1)})=N_1N_2$.

{\it Case III :}  $H_1=\{ e_1 \}$. Therefore $p_2$ is an isomorphism and there exist two possibilities: either $H_2=\{ e_2 \}$ and $p_1$ is also an isomorphism or $H_2\neq \{ e_2 \}$. 
If $p_1$ and $p_2$ are isomorphisms, then $N_1=N_2=N$ and $G=\Delta \subset \mathfrak S_{N}\times \mathfrak S_{N}$ up to automorphism of one of factors, and since $((1,\,2), (1,\,2)) \in G$, this automorphism must be inner. If $H_2\neq \{ e_2 \}$, then $p_1\circ p_2^{-1}: \mathfrak S_{N_2}\to \mathfrak S_{N_1}$ is an epimorphism (not isomorphism). Since $N_i>2$, $\mathfrak S_{N_2}$ must coincide with $\mathfrak S_{4}$, $\mathfrak S_{N_1}=\mathfrak S_{3}$ and $H_2=K_4$. The rest of the proof of this case will be left to the reader.

{\it Case IV :}  $N_1=4$ and $H_1=K_4$. The case $H_2=\mathfrak S_{N_2}$ is impossible. In fact, if we consider $p_1$, then we obtain $\mid G \mid = 4!N_2!$. On the other hand, if to consider $p_2$, then $\mid G \mid = 4N_2!$, a contradiction. 

The case $H_2=\mathfrak A_{N_2}$ is also impossible. In fact, if we consider $p_1$, then $\mid G \mid = 4!N_2!/2$. On the other hand, if to consider $p_2$, then $\mid G \mid = 4N_2!$, a contradiction. 

The case $H_2=\{ e_2 \} $ coincides (up to indexing) with one in Case III. 

The case $N_2=4$ and $H_2=K_4$ will be left to the reader.

To complete the proof of Proposition 2 we note that $\deg f_{1,2}=N_1N_2$ and there exists an irreducible component $X_{(1,1)}$ of $X$ such that $\deg f_{\mid X_{(1,1)}}=(G:St_{(1,1)})$. Hence, by Lemma 9, $\deg f_{\mid X_{(1,1)}}=N_1N_2$ and, consequently, $X$ is irreducible always except the case when $N_1=N_2=N$ and $G\simeq \Delta \subset \mathfrak S_{N}\times \mathfrak S_{N}$. But the exceptional case corresponds to one when $(S_1,f_1)$ and $(S_2,f_2)$ are equivalent. 
\begin{pred} 
$X$ is non-singular.
\end{pred}
\dvo . We need to check the smoothness of $X$ only at $z\in f_{1,2}^{-1}(B)$ such that $p_1=g_1(z)\in R_1\subset S_1$ and $p_2=g_2(z)\in R_2\subset S_2$. Put $f_{1,2}(z)=s$ and choose a small neighborhood $V_1\subset S_1$ of $p_1$ (resp. $V_2\subset S_2$ of $p_2$) and a small neighborhood $U\subset \mathbb P^2$ of $s$ such that $f_i(V_i)=U$ and, in the chosen neighborhoods, there exist local holomorphic coordinates for which equations defining $f_i$ have the simplest form. 

Let $s\in B$ be a non-singular point of $B$ or a node. Then 
$f_i :V_i\to U$ is given by 
 $$
\begin{array}{lll}
u_{i,1}^2 & =  & v_1\, ;\\
u_{i,2} & =  & v_2\, ,
\end{array}
$$
where $v_1=0$ is an equation of $B\cap U$ (or one of the branches of $B$ if $s$ is a node). 
Therefore $V_1\times _{U} V_2$ in  $V_1\times V_2$ is given by equations 
$$
\begin{array}{lll}
u_{1,1}^2 & =  & u_{2,1}^2\, ;\\
u_{1,2} & =  & u_{2,2}\, ,
\end{array}
$$
or, equivalently,  
$$
\begin{array}{lll}
u_{1,1} & =  & \pm u_{2,1}\, ;\\
u_{1,2} & =  & u_{2,2}\, .
\end{array}
$$
Hence  $V_1\times _{U} V_2$ consists of two irreducible non-singular components one of which corresponds to the sign $+$ and the other does to $-$. Therefore the normalization $\widetilde{V_1\times _{U} V_2}$ of $V_1\times _{U} V_2$ is the disjoint union of two non-singular surfaces. 

Let $s=s_j \in B$ be a cusp and $(x,y)$ local coordinates around $s$ chosen in section 1.3. Let $V_1\subset S_1$ (resp. for $ S_2$) be a neighborhood of $p_1=p_{1,j}$ such that $f_1=f_{1\mid V_1}: V_1\to U$ is a three-sheeted covering ramified along $R_1\cap V_1$. Put $Y=V_1\times _U V_2$ and let
$\widetilde Y$ be the normalization of $Y$. Denote by $g_i:\widetilde Y\to V_i$ and $f_{1,2}:\widetilde Y\to U$ the corresponding natural morphisms. Since $(V_1,f_1)$ and $(V_2,f_2)$ are equivalent, $f_{1,2}$ corresponds to the homomorphism 
$$\varphi _{1,2}=\varphi _1\times \varphi _2 : \pi _1 (U\setminus B)\to \Delta \subset \mathfrak S_{3}\times \mathfrak S_{3}\subset \mathfrak S_{9},$$ 
thus, $\varphi _{1,2}(\pi _1 (U\setminus B))$ acts on $\overline{9}\simeq \overline{3}\times \overline{3}$. It is easy to check that there are two orbits of the action of $\varphi _{1,2}(\pi _1 (U\setminus B))$: one of them is the orbit of $(1,1)$, and the other one is the orbit of $(1,2)$. Therefore $\widetilde Y$ is the disjoint union of $\widetilde Y_{(1,1)}$ and $\widetilde Y_{(1,2)}$. It is easily seen (cf. Lemma 1.6 in \cite{Cat}) that $(\widetilde Y_{(1,1)},f_{1,2})$ is isomorphic to $(\overline{V},\overline{f})$ (in notation of 1.3) and $g_i:\widetilde Y_{(1,1)}\to V_i$ is an isomorphism for $i=1,\, 2$;
$(\widetilde Y_{(1,2)},f_{1,2})$ is isomorphic to $(\overline{W},\widetilde{f})$, and each $(\widetilde Y_{(1,2)}, g_i)$ is isomorphic to $(\overline{W},g)$. Hence $X$ is non-singular.

\begin{zam} If $(S_1,f_1)$ and $(S_2,f_2)$ with the same discriminant curve are equivalent, then 

1) $X$ is non-singular;

2) $X$ is the disjoint union of two irreducible components: $X= X_{(1,1)}\bigsqcup X_{(1,2)}$, such that $g_{i|X_{(1,1)}}:X_{(1,1)}\to S_i$ is an isomorphism, $i=1,\, 2$, and $\deg g_{i|X_{(1,2)}}=N-1$, where $N=N_1=N_2$. 
\end{zam}
{\bf 2.2.} Let $\widetilde R \subset X$ be a curve $g_1^{-1}(R_1)\cap g_2^{-1}(R_2)$, $\widetilde C=g_1^{-1}(C_1)\cap g_2^{-1}(C_2)$, $\widetilde C_1=g_1^{-1}(R_1)\cap g_2^{-1}(C_2)$, and $\widetilde C_2=g_1^{-1}(C_1)\cap g_2^{-1}(R_2)$.
\begin{pred}
$$
\begin{array}{cll}
\widetilde R^2 & = & 2(3d+g-1)-c \, ,  \\
\widetilde C_1^2 & = & (N_2-2)(3d+g-1)-c   \, , \\
\widetilde C_2^2 & = & (N_1-2)(3d+g-1)-c   \, , \\
(\widetilde R,\widetilde C_i) & = & c\, \hspace{1cm} \mbox{for}\, \, i=1,\, 2 .
\end{array}
$$
\end{pred}
\dvo . It follows from the local considerations in the proof of Proposition 3 that $\deg g_{1\mid \widetilde R}=2$ and $g_{1\mid \widetilde R}$ is \'{e}tale.

One can check that $\widetilde R$ and $\widetilde C_i$
are intersected only at points over the cusps of $B$. Consider one of cusps of $B$, say $s=s_j$, and $p_i=p_{i,j}\in R_i\cap f^{-1}(s)$. In notation of the proof of Proposition 3, let $U$ be a neighborhood of $s$. It is easily seen that one of the branches of $\widetilde R\cap \widetilde Y$ lies in 
$\widetilde Y_{(1,1)}$ and the other one does in $\widetilde Y_{(1,2)}$. Since $(\widetilde Y_{(1,1)},f_{1,2})$ is isomorphic to $(V_i,f_i)$, the intersection $\widetilde Y_{(1,1)}\cap \widetilde C_i=\emptyset $ for $i=1,\, 2$. Consider  $\widetilde Y_{(1,2)}\cap \widetilde R$ and $\widetilde Y_{(1,2)}\cap \widetilde C_i$. Since each $(\widetilde Y_{(1,2)}, g_i)$ is equivalent to $(\overline{W},g)$, where $(\overline{W},g)$ was defined in section 1.3, we can identify $(\widetilde Y_{(1,2)}, g_1)$ with $(\overline{W},g)$. Then $\widetilde Y_{(1,2)}\cap \widetilde R$, $\widetilde Y_{(1,2)}\cap \widetilde C_1$, and $\widetilde Y_{(1,2)}\cap \widetilde C_2$ can be identified, respectively, with $\widetilde R$,  $\widetilde C_1$, and  $\widetilde C_2\subset \overline{W}$. Since in the neighborhood $\overline{W}$ the intersection multiplicity of $\widetilde R$ and $\widetilde C_i$ is equal to 1, hence $(\widetilde R,\widetilde C_i)=c$. 

To calculate $\widetilde R^2$, consider again the local case. Let $g:\overline W\to V$ be a two-sheeted covering given in local analytic coordinates by
 $$
\begin{array}{lll}
w_{1}^2 & =  & v_1\, ;\\
w_{2} & =  & v_2\, ,
\end{array}
$$
Denote by $C\subset V$ a curve given by $v_1=0$ and $R$ given by $v^2_2=v_1$. Then $C$ is the branch curve, $g^{*}(C)=2\widetilde C_2$, and $g^{*}(R)=\widetilde R+\widetilde C_1$, where $\widetilde C_2$ is given by $w_1=0$, and $\widetilde R$ and $\widetilde C_1$ are given by $w_2=\pm w_1$. Let $\sigma : \widetilde V\to V$ be the composition of two $\sigma $-processes with centers at points and such that $\sigma ^{-1}(R+C)=R+C+L_1+L_2$ is a divisor with normal crossings, where $L_1$ is the exceptional divisor of the first $\sigma $-process, $L_2$ is the exceptional divisor of the second $\sigma $-process, and, for simplicity of notation, we again denote by $R$ and $C$, respectively, the strict preimages of $R\subset V$ and $C\subset V$. Since we performed two $\sigma $-processes with centers at points lying in $R$, $R^2$ is decreased by 2 (if $R$ is considered as a complete curve). We can perform two $\sigma $-processes $\widetilde \sigma :\widetilde W \to \overline W$ (the fi!
rst one with center at the origin and the second one with center at the intersection point of the strict preimage of the curve $\{ w_2=0 \} $ and the exceptional divisor of the first $\sigma $-process). It is easy to check that  we again obtain a morphism $\overline g:\widetilde W\to \widetilde V$. Since we performed only one $\sigma $-process with center at a point lying in $\widetilde R$, $\widetilde R^2$ is decreased by 1. Besides, $\overline g_{\mid \widetilde R} :\widetilde R\to R$ is an isomorphism (locally), and $\overline g$ is non-ramified at each point lying in $\widetilde R$.

The considerations described above allow us to calculate $\widetilde R^2$. Indeed, performing at each point $p_{1,j}\in R_1$ two $\sigma $-processes as above, $R^2_1$ is decreased by $2c$.
Performing at each point of $g_1^{-1}(p_{1,j})\cap \widetilde R$ two $\sigma $-processes either as above or if a neighborhood of the point in consideration is isomorphic to $\widetilde Y_{(1,1)}$ we perform $\sigma $-processes as in $V_1$ as, in view of the fact that $g_1:\widetilde Y_{(1,1)}\to V_1$ is an isomorphism. After performing all these $\sigma $-processes, $\widetilde R^2$ is decreased by $3c$, and 
we can find a neighborhood $\widetilde V_1$ of the strict preimage of $R_1$ and a neighborhood $\widetilde W$ of the strict preimage of $\widetilde R$ such that the restriction $\widetilde g_{1\mid \widetilde W}:\widetilde W\to \widetilde V$ of the obtained morphism $\widetilde g_1$ to $\widetilde W$ is a non-ramified two-sheeted covering. Hence 
$$\widetilde R^2-3c=2(R_1^2-2c).$$  
Thus, by Lemma 4, $\widetilde R^2=2(3d+g-1)-c$.

Since $\deg g_1 =N_2$, 
$$N_2R_1^2=(\widetilde R +\widetilde C_1,\widetilde R+\widetilde C_1)= \widetilde R^2 +2 (\widetilde R,\widetilde C_1)+\widetilde C_1^2.$$
Hence
$$\widetilde C_1^2 =(N_2-2)(3d+g-1)-c .$$ Proposition 4 is proved. \newline 
{\bf 2.3.} To complete the proof of Theorem 1, we apply Hodge's Index Theorem. Since by Corollary 1, $\widetilde R^2 > 0$,
$$
\left|
\begin{array}{cc}
\widetilde R^2 & (\widetilde R,\widetilde C_i) \\
(\widetilde C_i,\widetilde R) &   \widetilde C_i^2  
\end{array}
\right| 
=\left|
\begin{array}{cc}
2(3d+g-1)-c & c \\
c &   (N_j-2)(3d+g-1)-c  
\end{array}
\right|\leq 0 ,
$$
i.e. 
$$2(N_j-2)(3d+g-1)^2-N_j(3d+g-1)c \leq 0 \, .$$
Hence 
$$N_j[2(3d+g-1)-c] \leq 4(3d+g-1) \, .$$ 
Thus, if there exist two non-equivalent generic morphisms $(S_1,f_1)$ and $(S_2,f_2)$, then their degrees
$$N_j \leq \frac{4(3d+g-1)}{2(3d+g-1)-c}.$$

\section{Uniqueness of generic morphisms for certain types of discriminant curves.}
Let us write inequality (\ref{in}) in the form 
\begin{equation}
N[2(3d+g-1)-c]-4(3d+g-1)>0 \label{ineq} 
\end{equation}
{\bf 3.1.} 
{\it Proof of Theorem} 2. Put $k=K_S^2$, then $E^2=m^2k$ and hence
\begin{equation}
N=\deg f=m^2k. \label{deg}
\end{equation}
We have $K_S=f^{*}(K_{\mathbb P^2})+R=-3E+R$, hence $R\equiv (3m+1)K_S$. Therefore $\deg B=(E,R)=m(3m+1)k$, i.e.
\begin{equation}
d=\frac{m(3m+1)k}{2} . \label{d}
\end{equation} 
By adjunction formula, $2(g-1)=(K_S+R,R)$. Therefore
\begin{equation}
g-1=\frac{(3m+2)(3m+1)k}{2} . \label{g}
\end{equation} 
{}From (\ref{d}) and (\ref{g}) it follows that
\begin{equation}
3d+g-1=(3m+1)^2k . \label{()}
\end{equation} 
By Lemma 8,
\begin{equation}
c=(12m^2+9m+3)k-12p_a . \label{c}
\end{equation}
Let us substitute (\ref{deg}), (\ref{()}), and (\ref{c}) into (\ref{ineq}). We have 
$$
\begin{array}{lll}
\displaystyle m^2k[(6m^2+3m-1)k +12p_a]-4(3m+1)^2k & = & \\
\displaystyle m^2k[(6m^2+3m-1)k +12p_a-4(3+\frac{1}{m})^2] & > & 0 .
\end{array}
$$
If $m\geq 2$, then 
$$
[(6m^2+3m-1)k +12p_a-4(3+\frac{1}{m})^2] \geq 29k+12p_a-49>0,
$$
and, since $k\geq 1$ and $p_a\geq 1$, the last inequality does not hold only for $k=1$ and $p_a=1$. In the exceptional case if $m\geq 3$, then
$$
[(6m^2+3m-1)k +12p_a-4(3+\frac{1}{m})^2] > 62+12-45>0.
$$
The case $m=2$, $k=p_a=1$ is impossible. Indeed, in this case by (\ref{d}) - (\ref{c}), we have
$\deg B=14$, $g=29$, and $c=57$, and the number of nodes $n$ must be non-negative. But by adjunction formula, 
$$n=\frac{1}{2}(\deg B-1)(\deg B-2)-g-c=13\cdot 6-29-57<0,$$
a contradiction. Therefore inequality (\ref{in}) holds always for $m\geq 2$.

If $m=1$, then 
$$
[(6m^2+3m-1)k +12p_a-4(3+\frac{1}{m})^2] = 8k+12p_a-64>0,
$$
and, consequently, inequality (\ref{in}) is equivalent to
$$2k+3p_a>16.$$
Since $k=E^2=N \geq 3$, the last inequality does not hold only for

1) $k=3, \, \, p_a\leq 3$;

2) $k=4, \, \, p_a\leq 2$;

3) $k=5, \, \, p_a\leq 2$;

4) $k=6, \, \, p_a=1$.

The cases 1) and 2) are impossible. Indeed, by (\ref{d}) - (\ref{c}), 
$\deg B=4k$, $g=10k+1$, and $c=12(2k-p_a)$, and the number of nodes $n$ must be non-negative. But by adjunction formula, 
$$n=\frac{1}{2}(\deg B-1)(\deg B-2)-g-c=(4k-1)(2k-1)-10k-1-24k+12p_a=4(2k^2-10k +3p_a)<0$$
for $k=3, \, \, p_a\leq 3$ and $k=4, \, \, p_a\leq 2$. 

The assertion on the degree of generic morphisms in the exceptional cases follows from Lemma 5.
\newline 
{\bf 3.2.} {\it Proof of Theorem} 3. Put $k=K_S^2$, then $E^2=m^2k$ and hence
\begin{equation}
N=\deg f=m^2k. \label{deg1}
\end{equation}
We have $K_S=f^{*}(K_{\mathbb P^2})+R=-3E+R$, hence $R\in |(-3m+1)K_S|$. Therefore $\deg B=(E,R)=m(3m-1)k$, i.e.
\begin{equation}
d=\frac{m(3m-1)k}{2} . \label{d1}
\end{equation} 
By adjunction formula, $2(g-1)=(K_S+R,R)$. Therefore
\begin{equation}
g-1=\frac{(3m-2)(3m-1)k}{2} . \label{g1}
\end{equation} 
It follows from (\ref{d1}) and (\ref{g1}) that
\begin{equation}
3d+g-1=(3m-1)^2k . \label{()1}
\end{equation} 
By Lemma 8,
\begin{equation}
c=(12m^2-9m+3)k-12 . \label{c1}
\end{equation}
Let us substitute (\ref{deg1}), (\ref{()1}), and (\ref{c1}) into (\ref{ineq}). We have 
$$
\begin{array}{lll}
\displaystyle m^2k[2(3m-1)^2k-(12m^2-9m+3)k +12]-4(3m-1)^2k & = & \\
\displaystyle m^2k[(6m^2-3m-1)k +12-4(3-\frac{1}{m})^2] &  = &  \\
\displaystyle m^2k[(6(m-\frac{1}{4})^2-\frac{11}{8})k +12-4(3-\frac{1}{m})^2] & > & 0 .
\end{array}
$$
If $m\geq 3$, then inequality (\ref{in}) holds, since
$$
[(6m^2-3m-1)k +12-4(3-\frac{1}{m})^2] > 26k+12-36>0.
$$
If $m=2$, then 
inequality (\ref{in}) holds, since
$$
[(6m^2-3m-1)k +12-4(3-\frac{1}{m})^2]=17k+12-25>0.
$$
If $m=1$, then 
inequality (\ref{in}) also holds, since in this case $k\geq 3$, hence
$$
[(6m^2-3m-1)k +12-4(3-\frac{1}{m})^2]=2k+12-16>0.
$$
{\bf 3.3.} {\it Proof of Theorem} 4. Let $f^{-1}(\mathbb P^1)=E\in |aL_1+bL_2|$, where $L_1$ and $L_2$ are the natural generators of $Pic\, S$. Without loss of generality, we can assume that $a\geq b>0$. Then
\begin{equation}
N=\deg f=2ab. \label{deg2}
\end{equation}
We have $K_S=-2L_1-2L_2$ and $K_S=f^{*}(K_{\mathbb P^2})+R=-3E+R$, hence $R\in |(3a-2)L_1+(3b-2)L_2|$. Therefore $\deg B=(E,R)=a(3b-2)+b(3a-2)$, i.e.
\begin{equation}
d=3ab-a-b . \label{d2}
\end{equation} 
By adjunction formula, $2(g-1)=(K_S+R,R)$. Therefore
\begin{equation}
g-1=9(ab-a-b)+8 . \label{g2}
\end{equation} 
It follows from (\ref{d2}) and (\ref{g2}) that
\begin{equation}
3d+g-1=18ab-12a-12b+8 . \label{()2}
\end{equation} 
By Lemma 7,
\begin{equation}
c=24ab-18a-18b+12 . \label{c2}
\end{equation}
Let us substitute (\ref{deg2}), (\ref{()2}), and (\ref{c2}) into (\ref{ineq}). We have 
$$
\begin{array}{lll}
\displaystyle 2ab[12ab-6a-6b+4]-4(18ab-12a-12b+8) & = & \\
\displaystyle 4ab[3a(b-1)+3b(a-1)+\frac{12}{b}+\frac{12}{a}-16-\frac{8}{ab}] & > & 0\, .
\end{array}
$$
If $a\geq b\geq 2$, then inequality (\ref{in}) holds, since
$$
3a(b-1)+3b(a-1)+\frac{12}{b}+\frac{12}{a}-16-\frac{8}{ab}\geq
3a+6(a-1)+\frac{8}{a}+\frac{12}{b}-16 =9a+\frac{8}{a}+\frac{12}{b}-22>0.$$
If $a>b=1$, then inequality (\ref{in}) also holds, since
$$
3a(b-1)+3b(a-1)+\frac{12}{b}+\frac{12}{a}-16-\frac{8}{ab}= 
3(a-1) +12+\frac{4}{a}-16 =3a+\frac{4}{a}-7>0.$$
If $a=b=1$, then $f$ is a two-sheeted covering of $\mathbb P^2$ branched along smooth conic. \newline
{\bf 3.4.} {\it Proof of Theorems} 5-7. Let 
\begin{equation}
\deg f=N= E^2=2k \label{deg3}
\end{equation}
 ($E^2$ is even, since $2K_S$ is trivial). We have $K_S=f^{*}(K_{\mathbb P^2})+R=-3E+R$, hence $R\equiv 3E$. Therefore $\deg B=6k$, i.e.
\begin{equation}
d=3k. \label{d3}
\end{equation} 
By adjunction formula, $2(g-1)=R^2$. Therefore
\begin{equation}
g-1=9k . \label{g3}
\end{equation} 
It follows from (\ref{d3}) and (\ref{g3}) that
\begin{equation}
3d+g-1=18k . \label{()3}
\end{equation} 
By Lemma 8,
\begin{equation}
c=24k-12p_a . \label{c3}
\end{equation}
Let us substitute (\ref{deg3}), (\ref{()3}), and (\ref{c3}) into (\ref{ineq}). We have 
$$
\begin{array}{lll}
2k[12k+12p_a]-72k & = & \\
24k[k+p_a-3] & > & 0\, .
\end{array}
$$
and $k\geq 2$, since $N>2$. If $S$ is a K3 surface, then inequality (\ref{in}) holds, since $p_a=2$.

If $S$ is an Enriques surface, then $p_a=1$ and inequality (\ref{in}) also holds except, possibly, for $k=2$. For abelian  varieties ($p_a=0$) inequality (\ref{in}) holds except, possibly, for $k=2$ and $k=3$. 

For abelian surfaces the case $k=2$ is impossible, since such a curve $B$ does not exist. Indeed, 
in this case it follows from (\ref{d3}), (\ref{g3}), and (\ref{c3}) that $\deg B=12$, $g=19$, $c=48$. Then $$n=\frac{1}{2}(\deg B-1)(\deg B-2)-g-c=55-19-48<0,$$ a contradiction. 

For Enriques surface in the exceptional case we have $\deg f=4$. Consider the epimorphism $\varphi : \pi _1 =\pi _1(\mathbb P^2\setminus B)\to \mathfrak{S}_4$ corresponding to $f$. It follows from (\ref{d3}), (\ref{g3}), and (\ref{c3}) that $\deg B=12$, $g=19$, and $c=36$, hence $n=0$. Therefore, by Proposition 1, the epimorphism $\varphi ^{\prime }: \pi _1 \to \mathfrak{S}_3$, which is the composition of $\varphi $ and 
the natural epimorphism $\mathfrak{S}_4\to \mathfrak{S}_3=\mathfrak{S}_4/K_4$, where $K_4$ is the Klein four group,
gives rise to a generic morphism $f^{\prime}: S^{\prime }\to \mathbb P^2$ of degree 3 with the same discriminant curve $B$. 

The assertion on the degree of generic morphisms in the exceptional cases follows from Lemma 5.
\newline
{\bf 3.5.} {\it Proof of Theorem} 8. Let $S=X(m_1,...,m_k)\subset \mathbb P^{k+2}$ be a complete intersection of multi-degree $(m_1,...,\, m_k)$, $m_i>1$. Then for a generic projection onto $\mathbb P^2$ 
\begin{equation}
\deg f=N= \prod _{i=1}^{k} m_i\, . \label{deg4}
\end{equation}
By adjunction formula, $K_S=(m_1+...+m_k-k-3)E$. Since $K_S=f^{*}(K_{\mathbb P^2})+R=-3E+R$, we have $R=(m_1+...+m_k-k)E$. Therefore $\deg B=(m_1+...+m_k-k)E^2$, i.e.
\begin{equation}
d=\frac{1}{2}(\sum _{i=1}^{k}(m_i-1))
\prod _{i=1}^{k} m_i\, . \label{d4}
\end{equation} 
By adjunction formula, $2(g-1)=R^2+(R,K_S)$. Therefore
\begin{equation}
g-1=\frac{1}{2}(\sum _{i=1}^{k}(m_i-1))(2\sum _{i=1}^{k}(m_i-1)-3)\prod _{i=1}^{k} m_i . \label{g4}
\end{equation} 
It follows from (\ref{d4}) and (\ref{g4}) that
\begin{equation}
3d+g-1= (\sum _{i=1}^{k}(m_i-1))^2
\prod _{i=1}^{k} m_i. \label{()4}
\end{equation} 
By Lemma 7,
\begin{equation}
c= 3N+2(g-1)-e(S)=3\prod _{i=1}^{k} m_i+(\sum _{i=1}^{k}(m_i-1))(2\sum _{i=1}^{k}(m_i-1)-3)\prod _{i=1}^{k} m_i-e(S). \label{c4}
\end{equation}
Let us substitute (\ref{deg4}), (\ref{()4}), and (\ref{c4}) into (\ref{ineq}). We have 
$$
\prod _{i=1}^{k} m_i[3\prod _{i=1}^{k} m_i(\sum _{i=1}^{k}(m_i-1))+e(S)-4(\sum _{i=1}^{k}(m_i-1))^2]>0 \, .
$$
If $k \geq 2$ , then 
$$
[3\prod _{i=1}^{k} m_i(\sum _{i=1}^{k}(m_i-1))+e(S)-4(\sum _{i=1}^{k}(m_i-1))^2]>(\sum _{i=1}^{k}(m_i-1))[3\prod _{i=1}^{k} m_i-4(\sum _{i=1}^{k}(m_i-1))]>0 
$$
and, obviously, the last inequality is true. 

If $k =1$ , then $e(S)=m_1^3-4m_1^2+6m_1$, hence 
$$
\begin{array}{lll}
\displaystyle 3\prod _{i=1}^{k} m_i(\sum _{i=1}^{k}(m_i-1))+e(S)-4(\sum _{i=1}^{k}(m_i-1))^2 & = & \\
3 m_1(m_1-1)+m_1^3-4m_1^2+6m_1-4(m_1-1)^2 & = & \\ m_1^3-5m_1^2+11m_1-4 & = & \\
(m_1-2)^3+(m_1-2)^2+3(m_1-2)+6 & \geq  & 6>0
\end{array}
$$
for $m\geq 2$, i.e. inequality (\ref{in}) also holds in the case $k=1$. 
\newline
{\bf 3.6.} 
{\it Proof of Theorem} 9.  Put $k=K_S^2$, $L^2=a$, and $(K_S,L)=b$. Note that $a>0$. If $E=mL$, then 
\begin{equation}
N=\deg f=m^2a. \label{deg9}
\end{equation}
We have $K_S=f^{*}(K_{\mathbb P^2})+R=-3E+R$, i.e. $R=K_S+3mL$. Therefore 
$$
\begin{array}{cll}
R^2 & = & 9m^2a+6mb+k, \\
(K_S,R) & = & 3mb+k.
\end{array}
$$
By adjunction formula, 
\begin{equation}
2(g-1)=(K_S+R,R)= 9m^2+9mb+2k  \label{g9}
\end{equation} 
and by Lemma 4, 
\begin{equation}
3d+g-1= 9m^2a+6mb+k. \label{()9}
\end{equation} 
By (\ref{deg9}), (\ref{g9}), and Lemma 7,
\begin{equation}
c=12m^2a+9mb+2k-e(S) . \label{c9}
\end{equation}
Let us substitute (\ref{deg9}), (\ref{()9}), and (\ref{c9}) into (\ref{ineq}). We have 
$$
\begin{array}{lll}
\displaystyle m^2a[6m^2a+3mb +e(S)]-4(9m^2a+6mb+k) & = & \\
\displaystyle m^2a[6m^2a+3mb +e(S)-36a-\frac{24b}{ma}-\frac{4k}{m^2a}] & > & 0 ,
\end{array}
$$
Now it is clear that there exist a constant $m_0$ such that the last inequality holds for $m\geq m_0$. 

The proof in the case $L=K_S$ coincides with one of Theorem 2.
\newline
{\bf 3.7.} 
{\it Proof of Theorem} 10. 
Let $B^{*}\subset \mathbb P^{2{*}}$ be a nodal curve of genus $g$, $\deg B^{*}=\delta $, and $B$ the dual curve of $B^{*}$. Then (cf. \cite{Cat}) $B$ is the discriminant curve of some generic morphism $f$ of degree $\delta $. In fact, let $S$ be the normalization of   
$$ X= \{ \, ((x_1,x_2,x_3),(y_1,y_2,y_3))\in \mathbb P^{2}\times \mathbb P^{2{*}}\, \, \mid \, \, \sum _{i}x_iy_i=0\, , \, \, (y_1,y_2,y_3))\in B^{*}\, \} $$
and $f$ be induced by the projection $pr_1: \mathbb P^{2}\times \mathbb P^{2{*}}\to \mathbb P^{2}$. It is clear that $(x_1,x_2,x_3)\notin B$ iff $\sum _{i}x_iy_i=0$ is not tangent to $B^{*}$, i.e. $f^{-1}((x_1,x_2,x_3))$ has exactly $\delta $ points and  
\begin{equation}
N=\deg f=\delta ,\label{deg5}
\end{equation}
hence $B$ is the branch curve of $f$ and it is clear that $f$ is a generic morphism.

It follows from Pl\"{u}cker's formulas that   
\begin{equation}
d=\delta+(g-1) \label{d5}
\end{equation} 
and 
\begin{equation}
c=3\delta+6(g-1). \label{c5}
\end{equation} 
Let us substitute (\ref{deg5}) - (\ref{c5}) into inequality (\ref{in})\, : 
$$
\delta >\frac{4(3\delta +4(g-1))}{3\delta +2(g-1)}=8-\frac{12\delta }{3\delta +2(g-1)} .
$$
Thus, inequality (\ref{in}) holds for $\delta \geq 8$. 

Consider the cases $\delta \leq 7$. 

The case $\delta =7$. Inequality (\ref{in}) does not hold iff $$
\frac{12\delta }{3\delta +2(g-1)} \leq 1 \Longrightarrow 9\delta \leq 2(g-1)\Longrightarrow g\geq 33 .
$$ 
On the other hand,
$$
g \leq \frac{(\delta-1)(\delta-2)}{2}=\frac{6\cdot 5}{2} < 33.
$$ 
Thus, in the case $\delta =7$ inequality (\ref{in}) is true.

The case $\delta =6$. Inequality (\ref{in}) does not hold iff $$
\frac{12\delta }{3\delta +2(g-1)} \leq 2 \Longrightarrow 3\delta \leq 2(g-1)\Longrightarrow g\geq 10 .
$$ 
On the other hand,
$$
g \leq \frac{(\delta-1)(\delta-2)}{2}=\frac{5\cdot 4}{2}=10.
$$ 
Thus, in the case $\delta =6$ inequality (\ref{in}) is not true  iff $\nu =0$, i.e., possibly, there exist two non-equivalent generic morphisms only when $B$  possesses the following invariants: $\deg B=30$, $g=10$, $c=72$, $n=324$. For such $B$, if there exist another (non-equivalent to $f$) generic morphism $f_1$ of $\deg f_1=N_1$, then
$$
N_1\leq \frac{4(3d +(g-1))}{6d +2(g-1)-c}=\frac{4(3\cdot 15+(10-1))}{6\cdot 15 +2(10-1)-72}=6.
$$
 
The case $\delta =5$. Inequality (\ref{in}) does not hold iff $$
\frac{12\delta }{3\delta +2(g-1)} \leq 3 \Longrightarrow \delta \leq 2(g-1)\Longrightarrow g\geq 4 .
$$ 
On the other hand,
$$
g \leq \frac{(\delta-1)(\delta-2)}{2}=\frac{5\cdot 4}{2}=6.
$$ 
Thus, in the case $\delta =5$ inequality (\ref{in}) is not true  iff either $\nu =0$, or $\nu =1$, or $\nu =2$, i.e., possibly, there exist two non-equivalent generic morphisms only in the following three cases: 

0) $\deg B=20$, $g=6$, $c=45$ 

1) $\deg B=18$, $g=5$, $c=39$. 

2) $\deg B=16$, $g=4$, $c=33$. \newline
In all cases 0) - 2) if there exist two non-equivalent generic morphisms $f$ and $f_1$ with the same discriminant curve $B$, then the computation, similar to one described above, gives rise to
$\deg f_1=N_1\leq 5$.

The case $\delta =4$. We have $g\leq 3$. If $g=3$, then $\deg B=12$ and $c=24$. It is easy to check that if there exist two non-equivalent generic morphisms $f$ and $f_1$ with such discriminant curve $B$, then $\deg f_1=N_1\leq 5$.  Let us show that such a curve can not be the discriminant curve of a generic morphism $f_1:S_1\to \mathbb P^2$ of degree 5. In fact, in this case by Lemmas 6 - 8, $K_{S_1}^2=-7$, $e(S_1)=-5$, and $p_a=-1$.  Hence $S_1$ is a ruled surface over a curve $C$ of genus $g=2$. Let $\overline{S}_1$ be a relatively minimal model of $S_1$. Then $e(\overline{S}_1)=-4$, hence $e(S_1)\geq-4$, a contradiction. 

The case $g\leq 2$ will be considered in the proof of Theorem 11.  \newline
{\bf 3.8.} {\it Proof of Theorem} 11. In the next subsection we shall prove Theorem 12. By that theorem, the Chisini Conjecture holds for a curve $B$ of $g\leq 1$, and if in the case $g=2$ there exist two non-equivalent generic morphisms $f_1$ and $f_2$, then $d\leq 3$. 

Consider the case $g=2$ and $d=3$. Then the inequality opposite to (\ref{in}) takes the following form
$$
\displaystyle N_i\leq \frac{40}{20-c}.
$$
Therefore, if $N_i\geq 5$, then
$$
\displaystyle 5\leq  \frac{40}{20-c} \Longleftrightarrow 12 \leq c\leq 19.
$$
On the other hand,
$$
c=(2d-1)(d-1)-g-n=8-n\leq 8,
$$
a contradiction. 

A curve $B$ with invariants $g=2$ and $\deg B=4$ (i.e. $d=2$) can not be a discriminant curve. Indeed, in this case either $c=1$, or $n=1$, which contradicts Corollary 2.

Consider the case $g=3$. By Theorem 12, if there exist two non-equivalent generic morphisms $f_1$ and $f_2$ such that $N_1\geq 5$, then $d\leq 6$. 
  By Lemmas 2 and 3,
\begin{equation}
\displaystyle \frac{3}{2}d+1\leq c\leq 3d +6,
\label{c and c}
\end{equation}
and the inequality opposite to (\ref{in}) takes the following form
\begin{equation}
\displaystyle 5\leq \frac{12d+8}{6d+4-c} \Longleftrightarrow  \frac{18d+12}{5}\leq c\leq 6d+3.
\label{www}
\end{equation}

If $d=6$, then by (\ref{c and c}), $10\leq c \leq 24$. On the other hand, by (\ref{www}), $c\geq 24$
and, consequently, $c=24$ and $n=28$. Thus, we obtain the case considered already in the proof of Theorem 10.

If $d=5$, then by (\ref{c and c}), $9\leq c \leq 21$. On the other hand, by (\ref{www}), $c\geq 21$
and, consequently, $c=21$ and $n=12$. By Pl\"{u}cker's formula, 
$\delta =10\cdot 9-2\cdot 12-3\cdot 21=3$, contrary to $g=3$. 

If $d=4$, then by (\ref{c and c}), $7\leq c \leq 18$. On the other hand, by (\ref{www}), $c\geq 17$
and, consequently, by Corollary 2,  $c=18$ and $n=0$. By Pl\"{u}cker's formula, 
$\delta =8\cdot 7-3\cdot 21\leq 0$, a contradiction.

If $d=3$, then by (\ref{c and c}), $6\leq c \leq 15$. On the other hand, by (\ref{www}), $c\geq 14$
and, consequently, by Corollary 2,  $c=15$, which contradicts the inequality $g=(2d-1)(d-1)-c-n\geq 0$.

If $d=2$, i.e. $\deg B=4$. Hence $B$ is non-singular. Thus, $B$ can not be a discriminant curve.\newline
{\bf 3.9.} {\it Proof of Theorem} 12. Consider again the inequality opposite to (\ref{in}) 
$$
N_i\leq \frac{4(3d +(g-1))}{6d +2(g-1)-c}.
$$
Thus, if there exist two non-equivalent generic morphisms such that one of them has degree $N_i\geq 5$, then
$$
\begin{array}{rll}
\displaystyle \frac{4(3d +(g-1))}{6d +2(g-1)-c} & \geq 5 & \Longrightarrow  \\
5c & \geq  18d+6(g-1)  & \Longrightarrow \,\, \mbox{(by Lemma 2) } \\
15d +15(g-1) & \geq 18d+6(g-1)  & \Longleftrightarrow \\
3(g-1) & \geq d . &  
\end {array}
$$

\section{Canonical discriminant curves.}
{\bf 4.1.} A curve $B$ is said to be {a(n $m$-)canonical discriminant curve} if $B$ is the discriminant curve of a generic morphism $f:S\to \mathbb P^2$ given by a linear subsystem $\{ E\} \subseteq |mK_S|$, $m\in \mathbb N$. 
 
Let $B\subset \mathbb P^2$ be a curve of even degree $2d$ with ordinary cusps and nodes as the only singularities and $\widetilde \nu :R\to B$ the normalization. Put $\mathfrak{e}=\widetilde \nu ^{-1}(\mathbb P^1\cap B)$, $\mathfrak{c} =2\sum '\widetilde \nu ^{-1}(s_i)$, and $\mathfrak{n}=\sum ''\widetilde \nu ^{-1}(s_i)$, where we denote by $\sum '$ (resp. by $\sum ''$) summation over all cusps (resp. nodes) $s_i\in B$.   
\begin{pred} \label{nec}
Let $B$ and $R$ be as above. If $B$ is a canonical discriminant curve, then 

($i$) \hspace{1cm} $\displaystyle 
\frac{2d}{g-3d-1}:=m\in \mathbb N \, , $ 

($ii$) \hspace{1cm} $\displaystyle 
\frac{(g-3d-1)^2}{3d+g-1}:=k\in \mathbb N \, , $

($iii$) \hspace{1cm} $\displaystyle 
\frac{4d^2}{3d+g-1}:=N\in \mathbb N \, , $

($iv$) \hspace{1cm} $\displaystyle 
N+\frac{3g-3-9d-c}{12}:=p_a\in \mathbb N \, , $

($v$) There exists a divisor $\mathfrak{k}\in Pic\, R$ such that $$ K_R=(3m+2)\mathfrak{k} \hspace{1cm} \mbox{and} \hspace{1cm} \mathfrak{e}=m\mathfrak{k}.$$
Besides, 
$$ \dim H^0(R,{\cal{O}}_R(r\mathfrak{k}))=\frac{r(r-1)}{2}k+p_a $$
for $r=2,\, ...\, , 3m$. 

($vi$) \hspace{1cm} $\displaystyle 
\mathfrak{c}+\mathfrak{n}=[(2d-6)m-2]\mathfrak{k} \, , $

\end{pred}
\dvo .  By (\ref{d}) and (\ref{g}),
$$
\left\{
\begin{array}{rll}
2g-2 & = & (9m^2+9m+2)k \\
6d & = & (9m^2 +3m)k 
\end{array} \right. \Longrightarrow 
$$
$$
\begin{array}{rlll}
g-1-3d & = & (3m+1)k & \Longrightarrow \\
m & = & \displaystyle \frac{2d}{g-3d-1} . & 
\end{array}
$$
Now ($ii$) follows from (\ref{d}), ($iii$) follows from (\ref{deg}), and ($iv$) does from (\ref{c})

The element $\mathfrak{k}\in Pic\, R$ is the restriction of $K_S$ to $R$. Since $K_S$ is ample, $$\dim H^1(S,{\cal{O}}_S(-rK_S))=0$$ 
for $r>0$. Therefore from the exact sequence 
$$0\longrightarrow {\cal{O}}_S((r-3m-1)K_S) \longrightarrow {\cal{O}}_S(rK_S) \longrightarrow {\cal{O}}_R(rK_S)\longrightarrow  0
$$ 
it follows that $H^0(S,{\cal{O}}_S(rK_S) )$ is isomorphic to $H^0(R,{\cal{O}}_R(rK_S))$ for $1\leq r\leq 3m$.  By Riemann-Roch formula, we have
$$\dim H^0(S,{\cal{O}}_S(rK_S))=\frac{r(r-1)}{2}K_S^2+p_a $$ for $r>1$. 

To prove ($vi$), let us blow up all singular points of $B$. Denote the composition of these $c+n\, \, \, \, 
\sigma$-processes by $\sigma : \widetilde{\mathbb{P}^2}\to \mathbb P^2$ and let $L_i=\sigma ^{-1}(s_i)$ be the exceptional curve over $s_i\in Sing \, B$. Then $\sigma ^{*}(B)=R+2\sum L_i$, where the strict preimage $R$ of $B$ is a non-singular curve, since all singularities of $B$ are ordinary cusps and nodes. We have 
$$ 
\begin{array}{rll}
K_{\widetilde{\mathbb{P}^2}} & = & -3\sigma ^{*}(\mathbb P^1)+\sum L_i\, , \\
R & = & \deg B\sigma ^{*}(\mathbb P^1)-2\sum L_i\, . 
\end{array}
$$
Hence by adjunction formula,
$${\cal{O}}_R(K_R)={\cal{O}}_R(K_{\widetilde{\mathbb{P}^2}}+R)=
{\cal{O}}_R((\deg B-3)\sigma ^{*}(\mathbb P^1)-\sum L_i).
$$
But ${\cal{O}}_R(\mathfrak{c}+\mathfrak{n})={\cal{O}}_R(\sum L_i)$. Therefore $${\cal{O}}_R(\mathfrak{c}+\mathfrak{n})={\cal{O}}_R(((2d-6)m-2)\mathfrak{k}),$$ 
since ${\cal{O}}_R(\sigma ^{*}(\mathbb P^1))={\cal{O}}_R(\mathfrak{e})={\cal{O}}_R(m\mathfrak{k})$ and ${\cal{O}}_R(K_R)={\cal{O}}_R((3m+2)\mathfrak{k})$ . 
\begin{pred} \label{plyu}
Let $B$ satisfy conditions ($i$) - ($vi$) of Proposition \ref{nec}. If $B$ is a discriminant curve of a generic morphism $f:S\to \mathbb P^2$ of $\deg f= N$, then $S$ is a surface of general type with ample canonical bundle and $f$ is given by a three-dimensional linear subsystem of $|mK_S+\alpha |$, where $\alpha \in Pic\, S$ is an element of order 2. Moreover, if $m$ is even, then $\alpha =0$.
\end{pred}
\dvo . Let $f$ be given by linear subsystem of $|E|$. We have
$$
E^2=N=\displaystyle \frac{4d^2}{3d+g-1} , 
$$ 
$R^2=3d+g-1$ and $(E,R)=2d$. Then 
$$
\left| \begin{array}{cc}
 E^2 &  (E,R) \\
(R,E) & R^2 
\end{array} \right| = 
\left| \begin{array}{cc}
 \displaystyle \frac{4d^2}{3d+g-1} & 2d \\
2d & {3d+g-1}  
\end{array} \right| =0.
$$
Therefore by Hodge's Index Theorem,  
$$R\equiv \displaystyle \frac{3d+g-1}{2d}E \equiv \frac{3m+1}{m}E.
$$
Since $K_S=-3E+R$, we have 
$$\begin{array}{rll}
E & \equiv  & mK_S,\\
R & \equiv & (3m+1)K_S,
\end{array}
$$ 
and if $E=mK_S+\alpha$, where $\alpha \equiv 0$, and 
${\cal{O}}_R(K_S)={\cal{O}}_R(\mathfrak{k}+\beta)$, where $\deg \beta = 0$, then $${\cal{O}}_R(E)={\cal{O}}_R(m\beta+\alpha+m\mathfrak{k})={\cal{O}}_R(m\mathfrak{k})$$
and 
$${\cal{O}}_R(R)={\cal{O}}_R(K_S+3E)={\cal{O}}_R((3m+1)\beta+3\alpha+(3m+1)\mathfrak{k}).$$
Hence ${\cal{O}}_R(K_R)={\cal{O}}_R((3m+2)\beta+3\alpha+(3m+2)
\mathfrak{k})$. Therefore
$$ \begin{array}{cll}
 {\cal{O}}_R(m\beta) & = & {\cal{O}}_R( -\alpha), \\
{\cal{O}}_R((3m+2)\beta) & = & {\cal{O}}_R( -3\alpha) . 
\end{array}  
$$
Hence $2\beta =0$ and ${\cal{O}}_R(2\alpha )={\cal{O}}_R$, and if $m$ is even, then ${\cal{O}}_R(\alpha )={\cal{O}}_R$.

Proposition \ref{plyu} will be proved once we prove  
\begin{lem} (cf. Appendix to Chapter V in \cite {zar2})
Let $S$ be a smooth projective surface and $i: R\hookrightarrow S$ a smooth irreducible curve. If ${\cal{O}}_S(R)$ is ample, then $i^{*}:Pic^0\, S\to Pic\, R$ is injective, where $Pic^0\, S\subset Pic\, S$ is the subgroup of numerically equivalent to zero classes of divisors.
\end{lem}
\dvo . In the commutative diagram 
$$
\xymatrix{
 &  &  0=H^1(S,{\cal{O}}_S(-R))\ar[d] & Pic^0\, S \ar@{^{(}->}[d] &  & \\
0 \ar[r] & H^1(S,\mathbb Z)\ar[d] \ar[r] & H^1(S,{\cal{O}}_S)\ar[r] \ar[d]^{j^{*}} & H^1(S,{\cal{O}}^{*}_S)\ar[r] \ar[d]^{i^{*}} & H^2(S,\mathbb Z)\ar[d] \ar[r]  & \\
0 \ar[r] & H^1(R,\mathbb Z) \ar[r] & H^1(R,{\cal{O}}_R)\ar[r]  & H^1(R,{\cal{O}}^{*}_R)\ar[r] & H^2(R,\mathbb Z) \ar[r] & }
$$
 with exact rows, the morphism $j^{*}$ is embedding. Therefore if for $\alpha \in Pic^0\, S$ the image 
$i^{*}(\alpha )=0$, then $r\alpha =0$ for some $r\in \mathbb N$, since $Tors\, H^2(S,\mathbb Z)$ is a finite abelian  group and $H^1(S,\mathbb Z)\to H^1(R,\mathbb Z)$ is embedding.  An element $\alpha \in  Pic\, S$ of finite order $r$ defines a non-ramified abelian  covering $\varphi: S_r\to S$ of $\deg \varphi =r$. If $\alpha \neq 0$ and  $i^{*}(\alpha)=0$, then $\varphi ^{-1}(R)$ is the disjoint union of $r$ irreducible curves, which contradicts the ampleness of $R$. 
\newline
{\bf 4.2.}  In notation of section 4.1, consider the natural homomorphism 
$$\mu : \mbox{Sym}^{2}H^0(R,{\cal{O}}_R(\mathfrak{e}))\longrightarrow H^0(R,{\cal{O}}_R(2\mathfrak{e})).$$ 
The kernel $\ker \mu$ generates the homogeneous ideal $I$ in the homogeneous coordinate ring $\mathfrak{R}= \oplus \mbox{Sym}^{r}H^0(R,{\cal{O}}_R(\mathfrak{e}))$ of the projective space 
$\mathbb P=\mathbb PH^0(R,{\cal{O}}_R(\mathfrak{e}))$. Put $S_I=Proj\, \mathfrak{R}/I$. The normalization $\widetilde \nu:R\to B$ determines a 3-dimensional subspace 
$$L=\widetilde \nu ^{*}(H^0(B,{\cal{O}}_B(\mathbb P^1\cap B)))\subseteq H^0(R,{\cal{O}}_R(\mathfrak{e})).$$ 
The subspace $L$ defines a projection $pr:\mathbb P \to \mathbb P^2$ with base locus $\mathbb PL\subset \mathbb P$. Let $f_I:S_I\to \mathbb P^2$ be the restriction of $pr$ to $S_I$.  
\begin{pred} \label{7}
Let $B$ satisfy conditions ($i$) - ($vi$) of Proposition \ref{nec}. Assume that $m$ is even and $\geq 21$. Then the curve 
$B$ is an $m$-canonical discriminant curve if and only if $S_I$ is a non-singular surface of $\deg S_I=N$ and $f_I$ is a generic morphism with discriminant curve $B$.
\end{pred}
\dvo . Let $A$ be a very ample divisor on $S$ and $D$ a numerically effective one. By \cite{Ein}, if the embedding of $S$ into $\mathbb P$ is given by $|K_S+4A+D|$, then the homogeneous ideal $I(S)$ is generated by quadrics. 

By \cite{Bom}, $A=5K_S$ is very ample (we assume that $K_S$ is ample). 
Therefore, if $S$ is embedded by $|mK_S|$, $m\geq 21$, then $I(S)$ is generated by quadrics. On the other hand, $R\in|(3m+1)K_S|$ and $R\subset S\subset \mathbb P$ is embedded by $|m\mathfrak{k}|$. As it was mentioned in the proof of Proposition \ref{nec}, the restriction map
$$H^0(S,{\cal{O}}_S(rmK_S))\longrightarrow H^0(R,{\cal{O}}_R(rmK_S))=H^0(R,{\cal{O}}_R(r\mathfrak{e}))$$ 
is an isomorphism for $r=1,\, 2$. Hence the set of quadrics containing $S$ coincides with the one containing $R$, and Proposition \ref{7} follows from Propositions \ref{nec} and \ref{plyu}.
\begin{zam} Proposition \ref{7} holds also for odd $m\geq 21$ if in view of Proposition \ref{plyu} we slightly change the definition of $m$-canonical discriminant curves. Moreover, Proposition 7 is a particular case of more general assertion.
\vskip 0.3cm 
\noindent {\bf Proposition $\bf{7}^{\prime }$} Let $B$ satisfy conditions ($i$) - ($iv$) of Proposition \ref{nec}. Assume that $m\geq 21$. Then the curve 
$B$ is the discriminant curve of morphism given by three-dimensional subsystem of $|L|$, where $L$ is numerically equivalent to the $m$th canonical class, if and only if $S_I$ is a non-singular surface of $\deg S_I=N$ and $f_I$ is a generic morphism with discriminant curve $B$.
\end{zam}
\section{On Zariski's pairs}
{\bf 5.1.} 
The set of plane curves of degree $2d$ is naturally parameterized by the points in $\mathbb P^{d(2d+3)}$. The subset of plane irreducible curves of degree $2d$ and genus $g$ with $c$ ordinary cusps and some nodes, as the only singularities, corresponds to a quasi-projective subvariety ${\cal{M}}(2d,g,c)\subset P^{d(2d+3)}$ (\cite{Wah}). One can show that if two non-singular points of the same irreducible component of ${\cal{M}}_{red}(2d,g,c)$ correspond to curves $B_1$ and $B_2$, then the pairs $(\mathbb P^2,B_1)$ and $(\mathbb P^2,B_2)$ are diffeomorphic. In particular, in this case the fundamental groups $\pi _1(\mathbb P^2\setminus B_1)$ and $\pi _1(\mathbb P^2\setminus B_1)$ are isomorphic.

The following Proposition is a simple consequence of Proposition 1 and local considerations in 1.2 and 1.3.

\begin{pred}
Let $(\mathbb P^2,B_1)$ and $(\mathbb P^2,B_2)$ be two diffeomorphic (resp. homeomorphic) pairs. If $B_1$ is the discriminant curve of a generic morphism $(S_1,f_1)$, then 
$B_2$ is also the discriminant curve of some generic morphism $(S_2,f_2)$. Moreover, if  $(S_1,f_1)$ is unique, then the same is true for $(S_2,f_2)$ and  $S_1$ and $S_2$ are diffeomorphic (resp. homeomorphic).
\end{pred}

Conversely, for $S\subset \mathbb P^r$ a projection $f:S\to \mathbb P^2$ is defined by a point in Grassmannian $\mbox{Gr}_{r+1,r-2}$ (the base locus of the projection). It is well-known that the set of generic projections is in one to one correspondence with some Zariski's open subset $U_S$ 
of $\mbox{Gr}_{r+1,r-2}$. A continuous variation of a point in $U_S$ gives rise to a continuous family of generic projections of $S$ whose branch curves belong to the same continuous family of plane cuspidal curves. Therefore discriminant curves of two generic projections of $S\subset \mathbb P^r$ belong to the same irreducible component of ${\cal{M}}(2d,g,c)$.
    
Moreover, if two surfaces $S_1$ and $S_2$ of general type with the same $K^2_S=k$ and $p_a=p$ are embedded by the $m$th canonical class into the same projective space $\mathbb P^r$ and belong to the same 
irreducible component of coarse moduli space ${\cal{M}}_S(k,p)$ of surfaces with given invariants (\cite{Gie}), then there exist generic projections $f_1$ of $S_1$ and $f_2$ of $S_2$ belonging to the same continuous family of generic projections. Therefore, discriminant curves of two generic projections of 
$S_1$ and $S_2$, belonging to the same 
irreducible component of a moduli space ${\cal{M}}_S(k,p)$, belong to the same irreducible component of ${\cal{M}}(2d,g,c)$ (cf. \cite{Wah}).
By Theorem 2 and Propositions 5 and 6, for a surface of general type with ample canonical class the triple of integers $(m,k,p)$ is uniquely determined by the invariants $(d,g,c)$ of $m$th canonical discriminant curve, and vice versa. Hence by Proposition 1 and Theorem 2, 
$$i(2d(k,p),g(k,p),c(k,p))\geq i(k,p),$$
where  $i(2d,g,c)$ (resp. $i(k,p)$) is the number of irreducible components of ${\cal{M}}(2d,g,c)$ (resp. ${\cal{M}}_S(k,p)$).
 
In \cite{Cat1}, F. Catanese showed that for each positive 
integer $h$ there exist integers $k,\, p$ such that ${\cal{M}}_S(k,p)$ has at least $h$ irreducible components. Hence for each positive integer $h$ there exist integers $d,\, g,\, c$ such that ${\cal{M}}(2d,g,c)$ has at least $h$ irreducible components. Note that the lower bound estimates
for $i(k,p)$, obtained in \cite {Cat1} and \cite{Man}, hold also for $i(2d,g,c)$, where $d,\, g$, and $c$ are the invariants of the corresponding $m$th canonical discriminant curves. 
\vskip 0.2cm 
\noindent 
{\bf 5.2.} 
Recently, there were published several articles (see, for example, \cite {Art}, \cite{Oka}, \cite{Shi}, \cite{Tok}) devoted to so called Zariski's pairs. By Artal-Bartolo's definition,  two plane curves $C_1,C_2\subset \mathbb P^2$ are called a Zariski pair if they have the same degree  and homeomorphic tubular neighborhoods in $\mathbb P^2$, but the pairs $(\mathbb P^2,C_1)$ and $(\mathbb P^2,C_2)$ are not homeomorphic. 

The first example of such pairs was obtained by O. Zariski and it is just two curves of degree 6 with 6 cusps mentioned in the Introduction.

In view of Theorem 2, applying Proposition 8, it is easy to prove the existence of a lot of Zariski pairs. To this end, we can find pairs of non-homeomorphic minimal surfaces of general type with ample canonical class and the same $K_S^2$ and $p_a$, and consider the corresponding $m$th canonical discriminant curves. For example,
\begin{pred}
For each integer $m\geq 5$, there is at least one Zariski's triple  (not only a pair) of plane cuspidal curves $B_{m,i}$, $i=3,\, 4,\, 5$, of degree $m(3m+1)$ and genus 
$g=\displaystyle \frac{1}{2}(3m+1)(3m+2)+1$ with $c=3(4m^2+3m-3)$ cusps.
\end{pred} 
\dvo . There exist (\cite{God}, \cite{Rei}, \cite{Bar}) five non-homeomorphic surfaces $S_i$ ($i=1,\, ...\, ,\, 5 $ ) of general type with $p_g=0$, $p_a=1$, and $K_{S_i}^2=1$. They are distinguished by $\mbox{Tors}H_1(S_i,\mathbb Z)=\mathbb Z/i\mathbb Z$. At least for $i\geq 3$ there exists a surface $S_{i}$ with ample canonical class possessing such invariants. Let $\phi _{m,i}:S_i\hookrightarrow \mathbb P^{\frac{1}{2}m(m-1)}$ be the $m$th canonical embedding, $m\geq 5$ (\cite{Bom}), and $f _{m,i}:S_i\to \mathbb P^2$ a generic projection with discriminant curve $B_{m,i}$. Applying the computation in the proof of Theorem 2, the discriminant curve $B_{m,i}$ has $$\deg B_{m,i}=m(3m+1), \hspace{0.5cm}
g=\displaystyle \frac{1}{2}(3m+1)(3m+2)+1,  \hspace{0.5cm} c=3(4m^2+3m-3).$$ By Proposition 8, since $S_i$ and $S_j$ are not homeomorphic for $i\neq j$, the pairs $(\mathbb P^2,B_{m,i})$ and $(\mathbb P^2,B_{m,j})$ also are not homeomorphic.

The easiest way to construct two non-homeomorphic surfaces $S_1$ and $S_2$ of general type, for which general morphisms, given by linear subsystems of the $m$th canonical class, give rise to Zariski's pair $B_1$ and $B_2$, is to construct two surfaces with different irregularities $q_1$ and $q_2$, where $q_i=\dim H^1(S_i,{\cal{O}}_{S_i})$, and the same $K_S^2$ and $p_a$. For example, let $A$ be an abelian surface and $C\subset A$ a non-singular curve such that the class of $C$ in $Pic\, A$ is divisible by 2. Consider the two-sheeted covering $\varphi :S\to A$ branched along $C$. It is easy to check that $S$ is a surface of general type with ample canonical bundle and irregularity $q\geq 2$. If $C^2=8p$, then the simplest computation gives rise to $K_S^2=4p$ and $p_a=p$. 
On the other hand, in \cite{Per}, U. Persson proved that for any 
positive integers $x$, $y$ satisfying 
\begin{equation}
2x-6\leq y\leq 8(x-cx^{2/3}),  \label{per}
\end{equation}
where $\displaystyle c=\frac{9}{\sqrt[3]{12}}$, there exists a simply connected minimal surface of general type with $K^2=y$ and $p_a=x$. It is easily seen that $x=p$ and $y=4p$ satisfy
 (\ref{per}) if $p\geq 486$.
\vskip 0.2cm 
\noindent  
{\bf 5.3.} {\bf Question}. {\it Let $S_1,S_2\subset \mathbb P^r$ be two diffeomorphic (resp. homeomorphic) surfaces of general type embedded by the $m$th canonical class. Are pairs $(\mathbb P^2,B_1)$ and $(\mathbb P^2,B_2)$ diffeomorphic (resp. homeomorphic), where $B_i$ is the discriminant curve of a generic projection of $S_i$ onto $\mathbb P^2$?}

Steklov Mathematical Institute

{\rm victor$@$olya.ips.ras.ru }

\end{document}